\documentclass[12pt]{article}%
\usepackage{amsfonts}
\usepackage{sw20bams}
\usepackage{amsmath}
\usepackage{amssymb}
\usepackage{graphicx}%
\providecommand{\U}[1]{\protect\rule{.1in}{.1in}}
\begin{document}

\title{Exercises in Iterational Asymptotics IV}
\author{Steven Finch}
\date{September 29, 2025}
\maketitle

\begin{abstract}
Abel's functional equation for $2^{x/2}$ and half-iterates of $\lambda
\,x(1-x)$ \& $\sqrt{1+x}$\ are featured in this collection of exercises
($0<\lambda\neq1<2$).

\end{abstract}

\footnotetext{Copyright \copyright \ 2025 by Steven R. Finch. All rights
reserved.}

This paper is a continuation of \cite{F1-exc4, F2-exc4, F3-exc4} and includes
work by\ Dmitrii Kouznetsov \& Henryk Trappmann (but employing an iterative
approach rather than power\ series).

\section{Onzi\`{e}me exercice}

\textbf{Consider }the recurrence
\[%
\begin{array}
[c]{ccccc}%
x_{k}=f(x_{k-1})=\sqrt{2}^{\;x_{k-1}} &  & \text{for }k\geq1\text{;} &  &
-\infty<x_{0}<\infty.
\end{array}
\]
Quantify the convergence rate of $x_{k}$ as $k\rightarrow\infty$ for both
$x_{0}=1$ and $x_{0}=3$. \ Let $F$ denote the solution, via $f$ iterates, of
Abel's equation $F(f(x))=F(x)+1$. \ Discover everything possible about $F$.
\ Using $F$, calculate $f^{[1/2]}(1)$, $f^{[1/2]}(3)$ and $f^{[1/2]}(5)$.

The function $f(x)=2^{x/2}$ has an attracting fixed point at $x=2$ and a
repelling fixed point at $x=4$; its inverse $g(x)=2\ln(x)/\ln(2)$ has a
repelling fixed point at $x=2$ and an attracting fixed point at $x=4$. \ For
$f$, set $y=x-2$ and note
\[%
\begin{array}
[c]{ccccc}%
y_{k}+2=2^{(y_{k-1}+2)/2} &  & \text{hence define} &  & \varphi(y)=2\left(
2^{y/2}-1\right)  ;
\end{array}
\]
for $g$, set $z=x-4$ and note%
\[%
\begin{array}
[c]{ccccc}%
z_{k}+4=\dfrac{2}{\ln(2)}\ln(z_{k-1}+4) &  & \text{hence define} &  &
\psi(z)=2\left(  \dfrac{1}{\ln(2)}\ln(z+4)-2\right)  .
\end{array}
\]
We have $\varphi^{\prime}(0)=\ln(2)$ and $\psi^{\prime}(0)=1/\ln(4)$, both of
which are between $0$ and $1$. \ The solution of Schr\"{o}der's equation
\[%
\begin{array}
[c]{ccc}%
\Phi\left(  \varphi(y)\right)  =\varphi^{\prime}(0)\Phi\left(  y\right)  , &
& \Psi\left(  \psi(z)\right)  =\psi^{\prime}(0)\Psi\left(  z\right)
\end{array}
\]
is obtained via Koenig's method \cite{Ax-exc4}:%
\[%
\begin{array}
[c]{ccc}%
\Phi\left(  y\right)  =\lim\limits_{k\rightarrow\infty}\dfrac{y_{k}}%
{\ln(2)^{k}}, &  & \Psi\left(  z\right)  =\lim\limits_{k\rightarrow\infty
}\dfrac{z_{k}}{\left(  1/\ln(4)\right)  ^{k}}%
\end{array}
\]
where $y=y_{0}$ and $z=z_{0}$. \ The function $\Phi$ helps to motivate%
\[
\lim\limits_{k\rightarrow\infty}\dfrac{x_{k}-2}{\ln(2)^{k}}=\left\{
\begin{array}
[c]{lll}%
-0.63209866105082925035545064599078086279947455... &  & \text{if }x_{0}=1,\\
2.18447475863901439313786713195265799616572364... &  & \text{if }x_{0}=3.
\end{array}
\right.
\]
The solution of Abel's equation%
\[%
\begin{array}
[c]{ccc}%
A\left(  \varphi(y)\right)  =A(y)+1, &  & B\left(  \psi(z)\right)  =B\left(
z\right)  +1
\end{array}
\]
is \cite{Ax-exc4}%
\[%
\begin{array}
[c]{ccc}%
A(y)=\dfrac{\ln\left(  \Phi(y)\right)  }{\ln(\varphi^{\prime}(0))}, &  &
B(z)=\dfrac{\ln\left(  \Psi(z)\right)  }{\ln(\psi^{\prime}(0))}%
\end{array}
\]
thus the desired function $F$ is%

\[
F(x)=\left\{
\begin{array}
[c]{lll}%
\dfrac{\ln\left(  -\Phi(x-2)\right)  }{\ln(\ln(2))}-\alpha &  & \text{if
}-\infty<x<2,\\
\dfrac{\ln\left(  \Phi(x-2)\right)  }{\ln(\ln(2))}-\beta &  & \text{if
}2<x<4,\\
\dfrac{\ln\left(  -\Psi(x-4)\right)  }{\ln(\ln(4))}-\gamma &  & \text{if
}2<x<4,\\
\dfrac{\ln\left(  \Psi(x-4)\right)  }{\ln(\ln(4))}-\delta &  & \text{if
}4<x<\infty
\end{array}
\right.
\]
where $\alpha$, $\beta$, $\gamma$, $\delta$ are constants (defined shortly).
\ To our surprise, there are two distinct solutions, supported on the interval
$(2,4)$, that are separated by an extremely short distance. \ The middle curve
in Figure 1 is thicker than the other curves, so as to convey that it
represents more than one solution. \ We discuss the four curves, moving from
left to right.

\subsection{$x$-intercept $1$}

Properties include%
\[%
\begin{array}
[c]{ccccccc}%
\lim\limits_{x\rightarrow-\infty}F(x)=-2^{+}, &  & F(0)=-1, &  & F(1)=0, &  &
\lim\limits_{x\rightarrow2^{-}}F(x)=\infty
\end{array}
\]
and%
\[
F(1)=\dfrac{\ln\left(  -\Phi(-1)\right)  }{\ln(\ln(2))}-\alpha=0
\]
forces%
\[
\alpha=1.25155147882218650957377135395164286460869893580054....
\]
This curve is identical to what was called $F_{2,1}^{-1}$ by Kouznetsov
\&\ Trappmann \cite{KT0-exc4}. \ Letting $T(c,0)=1$ and $T(c,n)=c^{T(c,n-1)}$,
we deduce that $F\left(  T\left(  \sqrt{2},n\right)  \right)  =n$. \ Although
a symmetry across the downward diagonal $y=-x$ is suggested by Figure 1, this
turns out to be false:\ intersections between the upward diagonal $y=x$ and
the curve are located at $x=-1.5459...$ and $x=1.5492...$.

\subsection{$x$-intercept $3_{a}$}

Properties include%
\[%
\begin{array}
[c]{ccccc}%
\lim\limits_{x\rightarrow2^{+}}F(x)=\infty, &  & F(3)=0, &  & \lim
\limits_{x\rightarrow4^{-}}F(x)=-\infty,
\end{array}
\]%
\[
F(5/2)=3.13739810096328698830281655519057645209387312849129...,
\]%
\[
F(7/2)=-3.23311619234714901269868926214873902220401094936075...
\]
and%
\[
F(3)=\dfrac{\ln\left(  \Phi(1)\right)  }{\ln(\ln(2))}-\beta=0
\]
forces%
\[
\beta=-2.13191778709502750839645694744655207295794772732641....
\]
This curve is identical to what was called $F_{2,3}^{-1}$ in \cite{KT0-exc4}.

\subsection{$x$-intercept $3_{b}$}

Properties include%
\[%
\begin{array}
[c]{ccccc}%
\lim\limits_{x\rightarrow2^{+}}F(x)=\infty, &  & F(3)=0, &  & \lim
\limits_{x\rightarrow4^{-}}F(x)=-\infty,
\end{array}
\]%
\[
F(5/2)=3.13739810096328698830281632291088083788688428250541...,
\]%
\[
F(7/2)=-3.23311619234714901269868829244855619555332830268389...
\]
and%
\[
F(3)=\dfrac{\ln\left(  -\Psi(-1)\right)  }{\ln(\ln(4))}-\gamma=0
\]
forces%
\[
\gamma=1.90057764535871549159209097160383434951546351559036....
\]
Values at $5/2$ and $7/2$ agree with the preceding to 24 decimal digits.
\ This curve is identical to what was called $F_{4,3}^{-1}$ in \cite{KT0-exc4}.

\subsection{$x$-intercept $5$}

Properties include%
\[%
\begin{array}
[c]{ccccc}%
\lim\limits_{x\rightarrow4^{+}}F(x)=-\infty, &  & F(5)=0, &  & \lim
\limits_{x\rightarrow\infty}F(x)=\infty
\end{array}
\]
and%
\[
F(5)=\dfrac{\ln\left(  \Psi(1)\right)  }{\ln(\ln(4))}-\delta=0
\]
forces%
\[
\delta=-1.11520724513160578643957308454984774926809594785349....
\]
This curve is identical to what was called $F_{4,5}^{-1}$ in \cite{KT0-exc4}. \ 

\subsection{Half-iterates}

The equation $F(x)=F(1)+1/2$ yields \cite{P-exc4}%
\[
f^{[1/2]}(1)=1.24362162766852180429509898360940293168819835661552...
\]
and the equation $F(x)=F(5)+1/2$ yields%
\[
f^{[1/2]}(5)=5.27364736917810646115785172418695262344663069400403....
\]
Two values of $f^{[1/2]}(3)$ are available \cite{AHS-exc4, K-exc4}:%
\[
f_{a}^{[1/2]}(3)=2.91378673463345260797944174759746005547139547790135...,
\]%
\[
f_{b}^{[1/2]}(3)=2.91378673463345260797944191697804065818551382121222...
\]
and these (as before) agree to 24 decimal digits. \ In terms of tetration,
$T\left(  \sqrt{2},1/2\right)  =1.2436...$.

We studied functions $e^{x/e}$ (with indifferent fixed point at $x=e$) and
$e^{x}$ (with no real fixed points at all) in \cite{F4-exc4, F5-exc4}, closely
following \cite{W-exc4, MS-exc4}. \ It was found that $T\left(  e^{1/e}%
,1/2\right)  =1.2571...$; using Kneser's construction, $T\left(  e,1/2\right)
=1.6463...$ but using Szekeres' construction, $T\left(  e,1/2\right)
=1.6451...$. \ Other examples via Kneser include $T\left(  3/2,1/2\right)
=1.2808...$, $T\left(  2,1/2\right)  =1.4587...$ and $T\left(  10,1/2\right)
=2.4770...$ \cite{PC-exc4, Ps-exc4}. \ It would be good to see the values of
$T\left(  4/3,1/2\right)  $, $T\left(  5/4,1/2\right)  $ and $T\left(
6/5,1/2\right)  $ someday.%

\begin{figure}
[ptb]
\begin{center}
\includegraphics[
height=6.0174in,
width=5.9871in
]%
{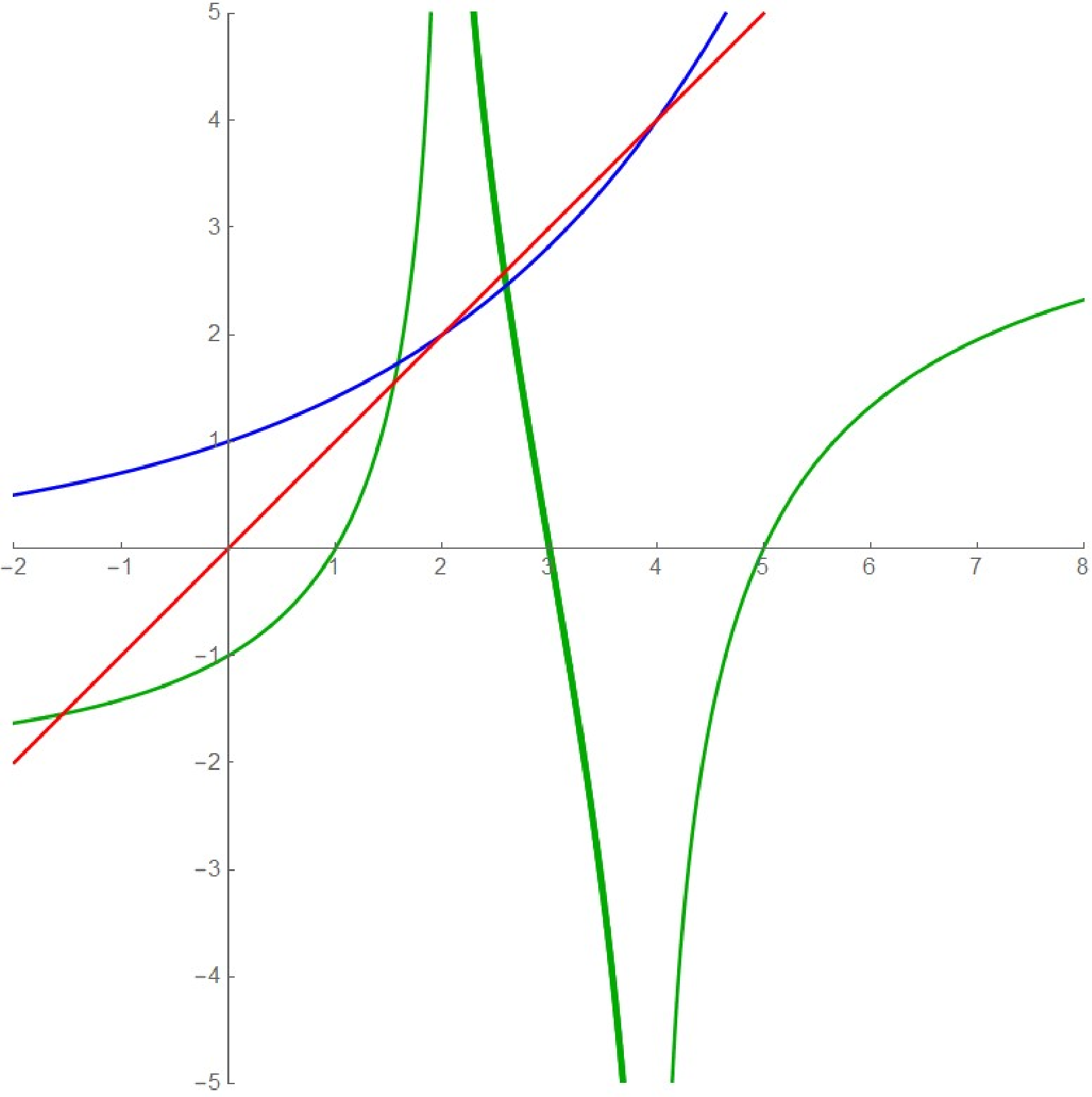}%
\caption{Four Abel solutions for $2^{x/2}$ in green (middle two are close,
hence thickness); $2^{x/2}$ in blue; upward diagonal in red.}%
\end{center}
\end{figure}
%

\begin{figure}
[ptb]
\begin{center}
\includegraphics[
height=5.9101in,
width=5.8738in
]%
{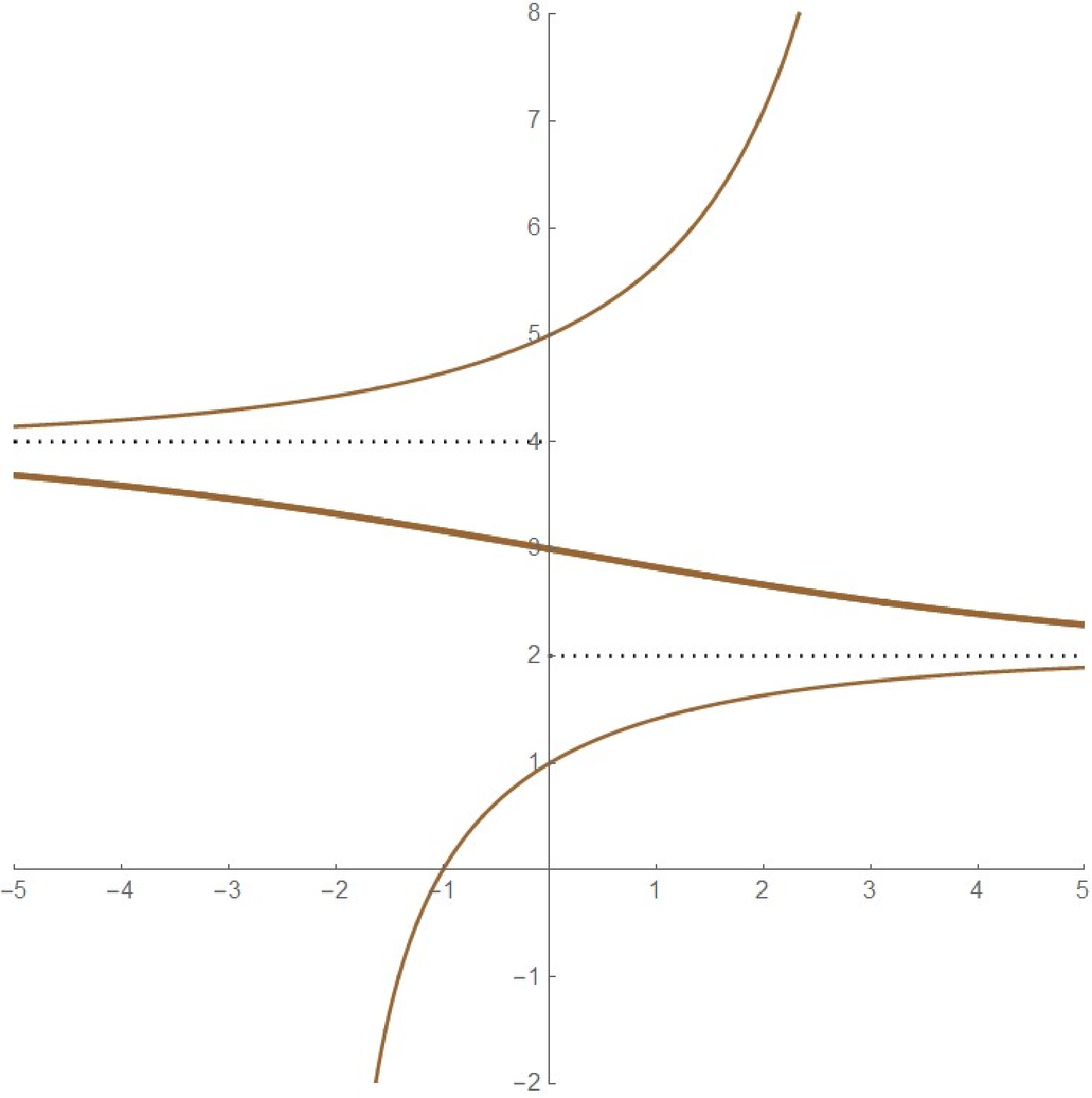}%
\caption{Inverses of Abel solutions, matching Figure 3 of Kouznetsov
\&\ Trappmann \cite{KT0-exc4} but omitting $10^{24}\,w(x)$ (where difference
$w(x)$ is approximately sinusoidal).}%
\end{center}
\end{figure}

\section{Douzi\`{e}me exercice}

\textbf{Consider initially }the recurrence \cite{F1-exc4}%
\[%
\begin{array}
[c]{ccccc}%
x_{k}=f(x_{k-1})=\lambda\,x_{k-1}(1-x_{k-1}) &  & \text{for }k\geq1\text{;} &
& x_{0}=1/2
\end{array}
\]
where $\lambda=1/2$, $1/3$ and $2/3$. \ Let $F$ denote the solution, via $f$
iterates, of Abel's equation $F(f(x))=F(x)+1$. \ Using $F$, calculate
$f^{[1/2]}(1/2)$ for each $\lambda$.

Solving this problem is straightforward. \ Because $0<f^{\prime}(0)=\lambda
<1$, Schr\"{o}der's equation $G(f(x))=f^{\prime}(0)G(x)$ gives%
\[%
\begin{array}
[c]{ccccc}%
G(x)=\lim\limits_{k\rightarrow\infty}\dfrac{x_{k}}{\lambda^{k}} &  &
\text{hence} &  & F(x)=\dfrac{\ln(G(x))}{\ln(\lambda)}.
\end{array}
\]
It follows that%
\[
f^{[1/2]}(1/2)=\left\{
\begin{array}
[c]{lll}%
0.2047412191193053167567775442679652620867325639... &  & \text{if }%
\lambda=1/2,\\
0.1628638902874469548727735573981779276762224789... &  & \text{if }%
\lambda=1/3,\\
0.2418762540161480676488939195406082513120285752... &  & \text{if }\lambda=2/3
\end{array}
\right.
\]
which are each less than $1/2$ and greater than $f(1/2)=1/8$, $1/12$, $1/6$
respectively. \ 

\textbf{Consider now} the same recurrence \cite{F1-exc4} where instead
$\lambda=3/2$, $4/3$ and $5/3$. \ As previously, calculate $f^{[1/2]}(1/2)$
for each $\lambda$.

The attracting fixed point here is $\mu=(\lambda-1)/\lambda$, not $0$. \ Set
$y=x-\mu$ and note%
\[%
\begin{array}
[c]{ccccc}%
y_{k}+\mu=\lambda(y_{k-1}+\mu)(1-y_{k-1}-\mu) &  & \text{hence define} &  &
g(y)=(2-\lambda)y-\lambda\,y^{2}.
\end{array}
\]
We have $0<g^{\prime}(0)=2-\lambda<1.$ \ Schr\"{o}der's equation
$G(g(y))=g^{\prime}(0)G(y)$ gives%
\[%
\begin{array}
[c]{ccccc}%
G(y)=\lim\limits_{k\rightarrow\infty}\dfrac{y_{k}}{(2-\lambda)^{k}} &  &
\text{thus} &  & F(x)=\dfrac{\ln(G(x-\mu))}{\ln(2-\lambda)}-\varepsilon
\end{array}
\]
where $\varepsilon$ is a constant depending on $\lambda$ (unlike earlier).
\ Assuming%
\[
F(\mu)=\dfrac{\ln(G(1/2-\mu))}{\ln(2-\lambda)}-\varepsilon=0
\]
forces%
\[
\varepsilon=\left\{
\begin{array}
[c]{lll}%
3.93270326530093996062495161372647676518537657293792... &  & \text{if }%
\lambda=3/2,\\
6.21293270406262080536684927084053547924538872086536... &  & \text{if }%
\lambda=4/3,\\
2.84080556393313600551740031635053919630961986852910... &  & \text{if }%
\lambda=5/3.
\end{array}
\right.
\]
It follows that%
\[
f^{[1/2]}(1/2)=\left\{
\begin{array}
[c]{lll}%
0.5984195936268982277477408185773449126377558120... &  & \text{if }%
\lambda=3/2,\\
0.6290618729919259661755530402296958743439857123... &  & \text{if }%
\lambda=4/3,\\
0.5674272219425106090254452885203644144647555042... &  & \text{if }\lambda=5/3
\end{array}
\right.
\]
which are each greater than $1/2$.

\textbf{Consider finally }the recurrence \cite{F6-exc4}%
\[%
\begin{array}
[c]{ccccc}%
x_{k}=f(x_{k-1})=\sqrt{1+x_{k-1}} &  & \text{for }k\geq1\text{;} &  & x_{0}=0
\end{array}
\]
i.e., the infinite simple radical expansion for the Golden mean $\varphi
=\left(  1+\sqrt{5}\right)  /2$. \ Let $F$ denote the solution, via $f$
iterates, of Abel's equation $F(f(x))=F(x)+1$. \ Using $F$, calculate
$f^{[1/2]}(0)$, $f^{[3/2]}(0)$ and $f^{[5/2]}(0)$.

The function $f(x)=\sqrt{1+x}$ has an attracting fixed point at $x=\varphi$.
\ Set $y=\varphi-x$ and note%
\[%
\begin{array}
[c]{ccccc}%
\varphi-y_{k}=\sqrt{1+\varphi-y_{k-1}} &  & \text{hence define} &  &
g(y)=\varphi-\sqrt{1+\varphi-y}.
\end{array}
\]
We have $0<g^{\prime}(0)=1/(2\varphi)<1.$ \ Schr\"{o}der's equation
$G(g(y))=g^{\prime}(0)G(y)$ gives%
\[%
\begin{array}
[c]{ccccc}%
G(y)=\lim\limits_{k\rightarrow\infty}\dfrac{y_{k}}{1/(2\varphi)^{k}} &  &
\text{therefore} &  & F(x)=\dfrac{\ln(G(\varphi-x))}{\ln(1/(2\varphi
))}-\varepsilon
\end{array}
\]
where $\varepsilon$ is a constant. \ Assuming
\[
F(0)=\dfrac{\ln(G(\varphi))}{\ln(1/(2\varphi))}-\varepsilon=0
\]
forces%
\[
\varepsilon=-0.67034187676403392725875840135990117450519933127639....
\]
The equation $F(x)=F(0)+n/2$ yields%
\[
f^{[n/2]}(0)=\left\{
\begin{array}
[c]{lll}%
0.58708229930179840752573065286737743155207668652537... &  & \text{if }n=1,\\
1.25979454646454094242043669146131942135336896191073... &  & \text{if }n=3,\\
1.50326130345477227760901754117049406417148121549859... &  & \text{if }n=5
\end{array}
\right.
\]
which are interspersed among $0$, $1$, $\sqrt{2}$ and $\sqrt{1+\sqrt{2}}$ as
expected. \ 

Analogous details corresponding to the simple continued fraction for $\left(
1+\sqrt{5}\right)  /2$ \&\ $1+\sqrt{2}$ \cite{F3-exc4}, the logistic with
$2<\lambda\leq3$ \cite{F2-exc4} and cosine \cite{F2-exc4} remain unresolved.
\ It is unclear whether oscillatory convergence (true for these open cases) is
compatible with Koenig's method.

\section{Addendum}

The two fixed points of $F$ in Section 1.1:%
\[
-1.54590582574454890961319276683302580776514319907836...,
\]%
\[
1.54921732984299390977708237603964130366954197864932...
\]
play a role for the function $\sqrt{2}^{\;x}$ comparable to what $\ell
\approx-1.85035452902718$ plays for the function $e^{x}$. \ The value $\ell$
is based on Kneser's construction; see formula (15)\ and Figure 5 of
\cite{Kz-exc4}. \ Another interesting constant $\approx1.63532$ is the largest
base $\tau$ (between $\sqrt{2}$ and $e$) in which the Abel solution for
$\tau^{x}$ has more than one fixed point.

A\ concluding example involves the factorial function \cite{KT1-exc4}:
$f(x)=\Gamma(1+x)$ has an attracting fixed point at $x=1$ and a repelling
fixed point at $x=2$; its inverse $g(x)$ for $x>0.8856$ has a repelling fixed
point at $x=1$ and an attracting fixed point at $x=2$.\ \ On the one hand, let
$z=x-2$ and note%
\[%
\begin{array}
[c]{ccccc}%
2+z_{k}=g(2+z_{k-1}) &  & \text{hence define} &  & h(z)=-2+g(2+z).
\end{array}
\]
We have $0<h^{\prime}(0)=1/(3-2\gamma)<1$, where $\gamma$ is Euler's constant.
\ Schr\"{o}der's equation $H(h(z))=h^{\prime}(0)H(z)$ gives%
\[%
\begin{array}
[c]{ccccc}%
H(z)=\lim\limits_{k\rightarrow\infty}\dfrac{z_{k}}{(1/(3-2\gamma))^{k}} &  &
\text{thus} &  & F(x)=\dfrac{\ln(H(x-2))}{\ln(3-2\gamma)}-\varepsilon
\end{array}
\]
where $\varepsilon$ is a constant. \ Assuming $F(3)=0$ forces%
\[
\varepsilon=-0.91938596545217952836264341194953987152246065907482...
\]
which improves upon a numerical estimate $-0.91938596545218$ on page 98 of
\cite{Kv-exc4}. \ Letting $w_{0}=3$ and $w_{n}=f(w_{n-1})$, we deduce that
$F(w_{n})=n$, i.e., $F$ increases without bound but extremely slowly. \ The
equation $F(x)=F(m)+1/2$ yields%
\[
f^{[1/2]}(3)=3.79606903179846431506883947400172995187762016400632...<6,
\]%
\[
f^{[1/2]}(4)=6.70253073232877914069364378501774180550013486665146...<24,
\]%
\[
f^{[1/2]}(5)=11.16011241010994435014306622367299550481488867202201...<120
\]
for integers $m=3$, $4$, $5$. \ On the other hand, let $y=x-1$ and note%
\[%
\begin{array}
[c]{ccccc}%
1+y_{k}=f(1+y_{k-1}) &  & \text{hence define} &  & g(y)=-1+f(1+y).
\end{array}
\]
We have $0<g^{\prime}(0)=1-\gamma<1.$ \ Schr\"{o}der's equation
$G(g(y))=g^{\prime}(0)G(y)$ gives%
\[%
\begin{array}
[c]{ccccc}%
G(y)=\lim\limits_{k\rightarrow\infty}\dfrac{y_{k}}{(1-\gamma)^{k}} &  &
\text{therefore} &  & F(x)=\dfrac{\ln(G(x-1))}{\ln(1-\gamma)}-\varepsilon.
\end{array}
\]
While $x=0$ is not a fixed point, it is mapped by $f$ to one. It is not
possible to determine $f^{[1/2]}(m)$ for any integer $-\infty<m\leq2$. \ We
turn attention to $x=-1/2$. \ Assuming $F(-1/2)=0$ forces%
\[
\varepsilon=-3.28913253388611161264763387262930308006425141579115...
\]
and the equation $F(x)=F(-1/2)+1/2$ yields%
\[
f^{[1/2]}(-1/2)=1.82377638928775824234619751192936424743765730719170...>\sqrt
{\pi}.
\]
Details about $F$ over the intervals $(m-1,m)$ for $0\neq m\leq2$ still linger.

\section{Acknowledgements}

The creators of Mathematica earn my gratitude every day:\ this paper could not
have otherwise been written. \ It constitutes perhaps my last words on the subject.


\begin{thebibliography}{99}                                                                                               %


\bibitem {F1-exc4}S. R. Finch, Exercises in iterational asymptotics, arXiv:2411.16062.

\bibitem {F2-exc4}S. R. Finch, Exercises in iterational asymptotics II, arXiv:2501.06065.

\bibitem {F3-exc4}S. R. Finch, Exercises in iterational asymptics III, arXiv:2503.13378.

\bibitem {Ax-exc4}D. S. Alexander, \textit{A History of Complex Dynamics. From
Schr\"{o}der to Fatou and Julia}, Friedr. Vieweg \& Sohn, 1994, pp. 46--49; MR1260930.

\bibitem {KT0-exc4}D. Kouznetsov and H. Trappmann, Portrait of the four
regular super-exponentials to base $\sqrt{2}$, \textit{Math. Comp.} 79 (2010)
1727--1756; MR2630010.

\bibitem {P-exc4}W. Paulsen, Tetration for complex bases, \textit{Adv. Comput.
Math.} 45 (2019) 243--267; MR3915009.

\bibitem {AHS-exc4}D. Asimov, D. Hickerson and R. Schroeppel, Are these two
numbers equal?, 1994 USENET math.sci newsgroup posting, http://groups.google.com/g/sci.math/c/HypjORccJqs.

\bibitem {K-exc4}J. Keane, Reply to \textquotedblleft Are these two numbers
equal?\textquotedblright, 1994 USENET math.sci newsgroup posting, http://groups.google.com/g/sci.math/c/HypjORccJqs.

\bibitem {F4-exc4}S. R. Finch, Half-iterates of $x(1+x)$, $\sin(x)$ and
$\exp(x/e)$, arXiv:2506.07625v1.

\bibitem {F5-exc4}S. R. Finch, Compositional square roots of $\exp(x)$ and
$1+x^{2}$, arXiv:2504.19999.

\bibitem {W-exc4}P. L. Walker, Infinitely differentiable generalized
logarithmic and exponential functions, \textit{Math. Comp.} 57 (1991)
723--733; MR1094963.

\bibitem {MS-exc4}K. W. Morris and G. Szekeres, Tables of the logarithm of
iteration of $e^{x}-1$, \textit{J. Austral. Math. Soc.} 2 (1961/62) 321--333; MR0141906.

\bibitem {PC-exc4}W. Paulsen and S. Cowgill, Solving $F(z+1)=b^{F(z)}$ in the
complex plane, \textit{Adv. Comput. Math.} 43 (2017) 1261--1282; MR3735895.

\bibitem {Ps-exc4}W. Paulsen, Tetration calculator, http://myweb.astate.edu/wpaulsen/tetcalc/tetcalc.html.

\bibitem {F6-exc4}S. R. Finch, Iterated radical expansions and convergence, arXiv:2410.02114.

\bibitem {Kz-exc4}D. Kouznetsov, Evaluation of holomorphic Ackermanns,
\textit{Appl. Comput. Math.} 3 (2014) 307--314.

\bibitem {KT1-exc4}D. Kouznetsov and H. Trappmann, Superfunctions and sqrt of
factorial, \textit{Moscow Univ. Physics Bull.} 65 (2010) 6--12.

\bibitem {Kv-exc4}D. Kouznetsov, \textit{Superfunctions}, 2020, https://mizugadro.mydns.jp/BOOK/466.pdf.%

\begin{tabular}
[c]{lll}
& Steven Finch & \\
& MIT Sloan School of Management & \\
& Cambridge, MA, USA & \\
& \textit{steven\_finch\_math@outlook.com} &
\end{tabular}

\end{thebibliography}
\end{document}